\documentclass[francais]{cedram-ambp}

\usepackage[autolanguage]{numprint}
\usepackage{amsmath,amsfonts,amssymb}
\usepackage{ragged2e}
\usepackage{natbib}
\usepackage[english,frenchb]{babel}
\usepackage[latin1]{inputenc}
\usepackage[T1]{fontenc}
 \usepackage{lmodern}
\usepackage{amsfonts}
\usepackage{amsmath}
\usepackage{pdfpages}

\usepackage{tikz}

\title
{La relation entre $\zeta(4n-1)$, $\zeta(2p)$ et $\zeta(4n-1-2p)$}

\author
[m. \lastname{kabongo}]
{\firstname{mundankulu} \lastname{kabongo}}

\address{Laboratoire de Recherche Lamunda\\
Sherbrooke/ Qu\'ebec / Canada}
\email{mundankulu.kabongo@gmail.com}
\begin{document}

%  Le r{\'e}sum{\'e}
\begin{abstract}
La relation fonctionnelle de la fonction z\^eta de Riemann ne nous fournit, ni la nature ni  l'expression, de z\^eta aux impairs positifs. A partir de la fonction $F(z)=\frac{z^{-2n}}{e^z-1}$ \`a variable complexe $z$, nous trouvons une relation fonctionnelle impliquant $\zeta(4n-1)$, $\zeta(2p)$ et $\zeta(4n-1-2p)$.Elle est donn\'ee par:

\begin{equation}
\zeta(4n-1)=\frac{1}{2n-1}\sum_{p=1}^{2n-2}\zeta(2p)\zeta(4n-1-2p).
\end{equation}\\$n=2, 3, 4, 5, 6, ...$\\ 
\\De cette formule nous introduisons une nouvelle approche pour \'etudier la nature de z\^eta \`a ces entiers.

\end{abstract}

% Le r�sum� en anglais
\begin{altabstract}
The functional relation of the Riemann z\^eta function provides us with neither the nature nor the expression of z\^eta at positive odd numbers. From the function $F(z)=\frac{z^{-2n}}{e^z-1}$, we find a functional relation involving $\zeta(4n- 1)$, $\zeta(2p)$ and $\zeta(4n-1-2p)$. It is given by:

\begin{equation}
\zeta(4n-1)=\frac{1}{2n-1}\sum_{p=1}^{2n-2}\zeta(2p)\zeta(4n-1-2p).
\end{equation}\\$n=2, 3, 4, 5, 6, ...$\\
\\From this formula we introduce a new approach to study the nature of z\^eta on these integers.
\end{altabstract}

%% Attention, les r�sum�s sont plac�s *avant* \maketitle
\maketitle 

 \newpage
%  D{\'e}but du texte de l'article
\section{Introduction}

En effet, la fonction z\^eta de Riemann ob\'eit \`a la relation fonctionnelle ${\pi}^{-\frac{s}{2}}\Gamma(\frac{s}{2})\zeta(s)={\pi}^{-\frac{1-s}{2}}\Gamma(\frac{1-s}{2})\zeta(1-s).$ Nous ne savons y tirer une expression ou la nature de  z\^eta aux entiers impairs positifs. Nous contournons la difficult\'e  qui apparait; en se servant de la fonction \`a variable complexe:
\begin{equation}
	F(z)=\frac{z^{-2n}}{e^z-1}
\end{equation}\\
Ci-bas, en d\'erivant $F(z)$, et au moyen d'une d\'emonstration, nous parvenons \`a obtenir la relation de $\zeta(4n-1)$; avec $n=2, 3, 4, ...$\\
\\ Puis nous avons introduit une nouvelle approche pour l'\'etude de la nature de $\zeta(3), \zeta(5),.,\zeta(4n-1)$ \\
Nous avions un demi cercle de rayon infini avec une infinit\'e des poles. Ces poles sont $z=2\pi j. p$ avec $p=1, 2, 3, 4, ...$. Nous avions trac\'e le demi cercle comme sur la figure ci-dessous, dans le but de pouvoir avoir une analyticit\'e de $F(z)$. Nous avons trouv\'e l'expression:

\begin{equation}
\zeta(4n-1)=\frac{1}{2n-1}\sum_{p=1}^{2n-2}\zeta(2p)\zeta(4n-1-2p)
\end{equation}\\$n=2, 3, 4, 5, 6, ...$\\
\\De cette relation,  nous obtenons:
\begin{equation}
 B_{4n-1}=\sum_{p=1}^{2n-2}\frac{A_{2p}}{2n-1}.B_{(4n-1-2p).}
\end{equation}\\Ceci nous permettra d'\'etudier z\^eta par les constantes $B_{4n-1}=\frac{\zeta(4n-1)}{{\pi}^{4n-1}}$, $\frac{A_{2p}}{2n-1}$ et $B_{4n-1-2p}=\frac{\zeta(4n-1-2p)}{{\pi}^{4n-1-2p}}$

\begin{equation}
A_{2p}=\frac{\zeta(2p)}{{\pi}^{2p}}
\end{equation}
 Ci-dessous la figure d'appui et la d\'emonstration de $\zeta(4n-1)$.
Notre domaine de calculs est $AabBGHIJA$. $F(z)$ est holomorphe dans ce domaine.\\ Le positif est compt\'e dans le sens antihoraire.
%\includepdf[pages=1]{domaine1}

\begin{tikzpicture}
%Axes
    \draw[black,thick,->] (-1,0) -- (12,0) node[anchor=north west] {X};
    \draw[black,thick,->] (5,-2) -- (5,12) node[anchor=north west] {Y};

% les barres des  liens entre les poles
    \draw[blue, line width=1mm, opacity=.8] (4.9,1.1) -- (4.9,1.9) node[anchor=north west] {};
    \draw[blue, line width=1mm, opacity=.8] (5.1,1.1) -- (5.1,1.9) node[anchor=north west] {};

  \draw[blue, line width=1mm, opacity=.5] (4.9,2.1) -- (4.9,2.9) node[anchor=north west] {};
    \draw[blue, line width=1mm, opacity=.5] (5.1,2.1) -- (5.1,2.9) node[anchor=north west] {};

  \draw[blue, line width=1mm, opacity=.5] (4.9,2.1) -- (4.9,3.9) node[anchor=north west] {};
    \draw[blue, line width=1mm, opacity=.5] (5.1,2.1) -- (5.1,3.9) node[anchor=north west] {};

  \draw[blue, line width=1mm, opacity=.5] (4.9,2.1) -- (4.9,4.9) node[anchor=north west] {};
    \draw[blue, line width=1mm, opacity=.5] (5.1,2.1) -- (5.1,4.9) node[anchor=north west] {};

  \draw[blue, line width=1mm, opacity=.8] (4.9,2.1) -- (4.9,6) node[anchor=north west] {};
    \draw[blue, line width=1mm, opacity=.8] (5.1,2.1) -- (5.1,6) node[anchor=north west] {};

% Les points  A  et B
\setlength{\unitlength}{1mm}
\put(-10,0){\circle*{2.5}}
\put(110,0){\circle*{2.5}}

\put(-15,-10){\makebox{$\Bigg(A\Bigg)$}}
\put(105,-10){\makebox{$\Bigg( B\Bigg)$}}

%Les points   a   et    b
\setlength{\unitlength}{1mm}
\put(48,0){\circle*{2}}
\put(52,0){\circle*{2}}

\put(43,-7){\makebox{$\Big(a\Big)$}}
\put(53,-7){\makebox{$\Big(b\Big)$}}

%Les points sur les barres des poles   J  et   G
\setlength{\unitlength}{1mm}
\put(49,60){\circle*{2}}
\put(51,60){\circle*{2}}

\put(40,63){\makebox{$\Big(J\Big)$}}
\put(53,63){\makebox{$\Big(G\Big)$}}

%Les points sur les barres des poles   I  et   H
\setlength{\unitlength}{1mm}
\put(49,35){\circle*{2.5}}
\put(51,25){\circle*{2.5}}

\put(40,37){\makebox{$\Big(I\Big)$}}
\put(53,27){\makebox{$\Big(H\Big)$}}

%Les points U   et    V
\setlength{\unitlength}{1mm}
\put(20,5){\circle*{2.5}}
\put(70,10){\circle*{2.5}}

\put(20,7){\makebox{$\Big(U\Big)$}}
\put(70,12){\makebox{$\Big(V\Big)$}}

%La partie du domaine qui est sur l'axe des X
 \draw[blue, line width=1mm, opacity=.5] (-1.5,0) -- (4.8,0) node[anchor=north west]{} ;
\draw[blue, line width=1mm, opacity=.5] (5.2,0) -- (12,0) node[anchor=north west] {};

%Le demi cercle
    \draw[blue, line width=1mm, opacity=.5] (5,0)+(0:6cm) arc[start angle=0, end angle=180, radius=6cm];

%Le petit demi cercle
      \draw[blue, line width=1mm, opacity=.8] (5,0)+(0:0.2cm) arc[start angle=0, end angle=180, radius=0.2cm];

%Les poles  
    \draw[blue, line width=1mm, opacity=.5] (5,1)+(0:0.2cm) arc[start angle=0, end angle=360, radius=0.2cm];
    \draw[blue, line width=1mm, opacity=.5] (5,2)+(0:0.2cm) arc[start angle=0, end angle=360, radius=0.2cm];
    \draw[blue, line width=1mm, opacity=.5] (5,3)+(0:0.2cm) arc[start angle=0, end angle=360, radius=0.2cm];
    \draw[blue, line width=1mm, opacity=.5] (5,4)+(0:0.2cm) arc[start angle=0, end angle=360, radius=0.2cm];
    \draw[blue, line width=1mm, opacity=.5] (5,5)+(0:0.2cm) arc[start angle=0, end angle=360, radius=0.2cm];

\end{tikzpicture}

\begin{remark}
Le demi cercle a un rayon \'egal \`a l'infini. Le point A a pour coordonn\'ees ($-\infty,0$) et le point B a pour coordonn\'ees ($+\infty,0$).
\end{remark}

\newpage
\section{D\'emonstration de $\zeta(4n-1)$ }
Prenons la fonction $F(z)$ et d\'erivons-la par rapport \`a $z$, on aura:
\\
\\
\begin{equation}
\frac{d}{dz}\Bigg(\frac{z^{-2n}}{e^z-1}\Bigg)=-2n\frac{z^{-2n-1}}{e^z-1}-\frac{z^{-2n}}{(e^z-1)^2}e^z
\end{equation}\\
%  Voici un exemple d'utilisation de l'environnement notation
%\begin{notation} 
 
%\end{notation}

%la figure que vous souhaitez insérer.

\justify
Comme $F(z)$ est holomorphe dans ce domaine, sa d\'eriv\'ee l'est aussi et on la note par petit $f$:\\
\\
\begin{equation}
f(z)=\frac{d}{dz}\Bigg(\frac{z^{-2n}}{e^z-1}\Bigg)=-2n\frac{z^{-2n-1}}{e^z-1}-\frac{z^{-2n}}{(e^z-1)^2}e^z
\end{equation}\\
\\Int\'egrons  sur le chemin $AUVB$, les deux membres de l'\'equation ci-dessus.On aura:

\begin{equation}
\int_{AUVB}f(z)dz=-2n\int_{AUVB}\frac{z^{-2n-1}}{e^z-1}dz-\int_{AUVB}\frac{z^{-2n}}{(e^z-1)^2}e^zdz
\end{equation}\\
\\Le point $A$  tend vers les coordonn\'ees $(-\infty, 0)$
\\Le point $B$  tend vers les coordonn\'ees $(+\infty, 0)$\\
\\
\\Le membre de gauche de (2.3) sera:

\begin{equation}
F(z)\Bigg|_{(-\infty,0)}^{(+\infty,0)}=0
\end{equation}\\
\\La relation (2.3) devient:

\begin{equation}
0=-2n\int_{AUVB}\frac{z^{-2n-1}}{e^z-1}dz-\int_{AUVB}\frac{z^{-2n}}{(e^z-1)^2}e^zdz
\end{equation}\\

\begin{equation}
2n\int_{AUVB}\frac{z^{-2n-1}}{e^z-1}dz=-\int_{AUVB}\frac{z^{-2n}}{(e^z-1)^2}e^zdz
\end{equation}\\
\\Puis on a:

\begin{equation}
\int_{AUVB}\frac{z^{-2n-1}}{e^z-1}dz=-\frac{1}{2n}\int_{AUVB}\frac{z^{-2n}}{(e^z-1)^2}e^zdz
\end{equation}\\
\\
\\On sait par ailleurs que:

\begin{equation}
Q_n=\oint_{AUVBGHIJA}\frac{z^{-2n-1}}{e^z-1}dz
\end{equation}\\  est nulle sur ce contour ferm\'e.\\
\\
\\Donc:
\begin{equation}
	0=\int_{AUVB}\frac{z^{-2n-1}}{e^z-1}dz+\int_{BG}\frac{z^{-2n-1}}{e^z-1}dz+\int_{GHIJ}\frac{z^{-2n-1}}{e^z-1}dz+\int_{JA}\frac{z^{-2n-1}}{e^z-1}dz
\end{equation}\\
\\Comme le rayon du demi cercle tend vers l'infini, alors le deuxi\`eme et le quatri\`eme terme sont nuls. Nous avons donc:
\begin{equation}
	0=\int_{AUVB}\frac{z^{-2n-1}}{e^z-1}dz+\int_{GHIJ}\frac{z^{-2n-1}}{e^z-1}dz
\end{equation}\\
\\
\\En r\'earrangeant la relation on a:
\begin{equation}
	\int_{AUVB}\frac{z^{-2n-1}}{e^z-1}dz=-\int_{GHIJ}\frac{z^{-2n-1}}{e^z-1}dz
\end{equation}\\En changeant le sens GHIJ en JIHG on a:

\begin{equation}
	\int_{AUVB}\frac{z^{-2n-1}}{e^z-1}dz=\int_{JIHG}\frac{z^{-2n-1}}{e^z-1}dz
\end{equation}\\
Le membre de droite de (2.12) est:

\begin{equation}
	(2\pi j)\sum_{p=1}^{\infty}(2\pi j. p)^{-2n-1}
\end{equation}\\
On a donc:
\begin{equation}
	\int_{AUVB}\frac{z^{-2n-1}}{e^z-1}dz=\frac{(-)^n}{(2\pi)^{2n}}.\zeta(2n+1)
\end{equation}\\
\\
\\En r\'eecrivant la relation (2.7) on a:
\begin{equation}
\int_{AUVB}\frac{z^{-2n-1}}{e^z-1}dz=-\frac{1}{2n}\int_{AUVB}\frac{z^{-2n}}{(e^z-1)^2}e^zdz
\end{equation}\\
La relation (2.14) dans (2.15) on a:
\begin{equation}
\frac{(-)^n}{(2\pi)^{2n}}.\zeta(2n+1)=-\frac{1}{2n}\int_{AUVB}\frac{z^{-2n}}{(e^z-1)^2}e^zdz
\end{equation}\\
\\
\\
Celui qui integre sur le chemin AUVB, c'est la m\^eme chose que celui qui integre surAabB, car cette int\'egrale ne d\'epend pas du chemin. Donc on peut \'ecrire:

\begin{equation}
\frac{(-)^n}{(2\pi)^{2n}}.\zeta(2n+1)=-\frac{1}{2n}\int_{AabB}\frac{z^{-2n}}{(e^z-1)^2}e^zdz
\end{equation}\\
Et on a:
$$
\int_{AabB}\frac{z^{-2n}}{(e^z-1)^2}e^zdz> 0
$$\\
Pour conserver l'\'egalit\'e des signes; \`a gauche et \`a droite de (2.17), alors il faut remplacer $n$ par $2n-1$

\begin{equation}
\frac{(-)^{2n-1}}{(2\pi)^{4n-2}}.\zeta(4n-1)=-\frac{1}{4n-2}\int_{AabB}\frac{z^{-4n+2}}{(e^z-1)^2}e^zdz
\end{equation}\\

\begin{equation}
\frac{1}{(2\pi)^{4n-2}}.\zeta(4n-1)=\frac{1}{4n-2}\int_{AabB}\frac{z^{-4n+2}}{(e^z-1)^2}e^zdz
\end{equation}\\

\begin{equation}
\frac{1}{(2\pi)^{4n-2}}.\zeta(4n-1)=\frac{1}{4n-2}\int_{AabB}\frac{z^{-4n+2}}{(e^z-1)}\Bigg(\frac{e^z}{e^z-1}\Bigg)dz
\end{equation}\\
\\
\\On sait que:

\begin{equation}
\Bigg(\frac{e^z}{e^z-1}\Bigg)=\sum_{m=0}^{\infty}e^{-m z}
\end{equation}\\
\\La relation (2.21) dans (2.20) alors on a:

\begin{equation}
\frac{1}{(2\pi)^{4n-2}}.\zeta(4n-1)=\frac{1}{4n-2}\int_{AabB}\frac{z^{-4n+2}}{(e^z-1)}\sum_{m=0}^{\infty}e^{-m z}dz
\end{equation}\\
\\Par ailleurs l'exponentielle $e^{-m z}$ s'\'ecrit:

\begin{equation}
e^{-m z}=\sum_{p=0}^{\infty}\frac{(-1)^p m^p z^p}{p!}
\end{equation}\\
\\La relation (2.23) dans (2.22) alors:

\begin{equation}
\frac{1}{(2\pi)^{4n-2}}.\zeta(4n-1)=\frac{1}{4n-2}\int_{AabB}\frac{z^{-4n+2}}{(e^z-1)}\sum_{m=0}^{\infty}\sum_{p=0}^{\infty}\frac{(-1)^p m^p z^p}{p!}dz
\end{equation}\\

On permute le signe int\'egrale avec les sommations; et on a:

\begin{equation}
\frac{1}{(2\pi)^{4n-2}}.\zeta(4n-1)=\frac{1}{4n-2}\sum_{m=0}^{\infty}\sum_{p=0}^{\infty}\frac{(-1)^p m^p }{p!}\int_{AabB}\frac{z^{-4n+2+p}}{(e^z-1)}dz
\end{equation}\\

Comme le membre de gauche de la relation ci-dessus est un nombre r\'eel, alors nous devons \'eliminer tous les complexes du membre de droite. Noter que les valeurs(termes) complexes proviennent des valeurs paires de la variable muette $p$. Nous \'ecrivons la relation seulement pour les $p$ impaires. On a donc:

\begin{equation}
\frac{1}{(2\pi)^{4n-2}}.\zeta(4n-1)=\frac{1}{4n-2}\sum_{m=0}^{\infty}\sum_{p=1}^{\infty}\frac{(-1)^{2p-1} m^{2p-1} }{(2p-1)!}\int_{AabB}\frac{z^{-4n+2+2p-1}}{(e^z-1)}dz
\end{equation}\\

Puis on \'elimine tous les couples $(m,p)$ tels que $m=0$, \`a cause que les termes correspondants sont nuls. La sommation des $m$ va commencer donc \`a $m=1$. On aura donc:

\begin{equation}
\frac{1}{(2\pi)^{4n-2}}.\zeta(4n-1)=\frac{1}{4n-2}\sum_{m=1}^{\infty}\sum_{p=1}^{\infty}\frac{(-1)^{2p-1} m^{2p-1} }{(2p-1)!}\int_{AabB}\frac{z^{-4n+1+2p}}{(e^z-1)}dz
\end{equation}\\

On \'ecrit donc:
\begin{equation}
\frac{1}{(2\pi)^{4n-2}}.\zeta(4n-1)=\frac{1}{4n-2}\sum_{p=1}^{\infty}-\frac{\zeta(-2p+1) }{(2p-1)!}\int_{AabB}\frac{z^{-4n+1+2p}}{(e^z-1)}dz
\end{equation}\\

On sait que:

\begin{equation}
\zeta(-2p+1)=\frac{-B_{2p}}{2p}
\end{equation}\\
(2.29) dans (2.28) alors on a:

\begin{equation}
\frac{1}{(2\pi)^{4n-2}}.\zeta(4n-1)=\frac{1}{4n-2}\sum_{p=1}^{\infty}\frac{B_{2p} }{(2p)!}\int_{AabB}\frac{z^{-4n+1+2p}}{(e^z-1)}dz
\end{equation}\\

L'int\'egrale dans cette relation devient:

\begin{equation}
\int_{AabB}\frac{z^{-4n+1+2p}}{(e^z-1)}dz=(2\pi j)^{-4n+2+2p}\zeta(4n-1-2p)
\end{equation}\\
\\L'int\'egrale ci-dessus a du sens seulement si $\alpha=-4n+1+2p < -1$. Ceci \`a cause de la conjonction de deux raisons:\\
\\
\\$\alpha=-1$ on a $\zeta(1)$(L'int\'egrale diverge)\\
\\
\\$\alpha>-1$ alors l'int\'egrale diverge  aussi, car le point $A$. est \`a moins l'infini.\\
\\
On  \'ecrit:$-4n+1+2p<-1$.\\
\\Cela implique $p\leq 2n-2$; on aura donc:
\\(2.31) dans (2.30) on a:

\begin{equation}
\frac{1}{(2\pi)^{4n-2}}.\zeta(4n-1)=\frac{1}{4n-2}\sum_{p=1}^{2n-2}\frac{B_{2p} }{(2p)!}(2\pi j)^{-4n+2+2p}\zeta(4n-1-2p)
\end{equation}\\

On a:
\begin{equation}
\frac{1}{(2\pi)^{4n-2}}.\zeta(4n-1)=\frac{1}{4n-2}\sum_{p=1}^{2n-2}\frac{B_{2p} }{(2p)!}(2\pi)^{-4n+2+2p}( j)^{-4n+2+2p}\zeta(4n-1-2p)
\end{equation}\\

\begin{equation}
\zeta(4n-1)=\frac{1}{4n-2}\sum_{p=1}^{2n-2}\frac{B_{2p} }{(2p)!}(2\pi)^{2p}( j)^{-4n+2+2p}\zeta(4n-1-2p)
\end{equation}\\

\begin{equation}
\zeta(4n-1)=\frac{1}{2n-1}\sum_{p=1}^{2n-2}\frac{B_{2p} }{2(2p)!}(2\pi)^{2p}( j)^{-4n+2+2p}\zeta(4n-1-2p)
\end{equation}\\

\begin{equation}
\zeta(4n-1)=\frac{1}{2n-1}\sum_{p=1}^{2n-2}\frac{B_{2p} }{2(2p)!}(2\pi)^{2p}( -1)^{-2n+1+p}\zeta(4n-1-2p)
\end{equation}\\

\begin{equation}
\zeta(4n-1)=\frac{1}{2n-1}\sum_{p=1}^{2n-2}\frac{B_{2p} }{2(2p)!}(2\pi)^{2p}( -1)^{1+p}\zeta(4n-1-2p)
\end{equation}\\

\begin{equation}
\zeta(4n-1)=\frac{1}{2n-1}\sum_{p=1}^{2n-2}\zeta(2p)\zeta(4n-1-2p)
\end{equation}\\

\section{Nouvelle approche pour \'etudier la nature de zeta aux impairs positifs}

D\'efinissons\\
\\
\begin{equation}
 \frac{\zeta(4n-1)}{{\pi}^{4n-1}}=B_{4n-1}
\end{equation}

\begin{equation}
 \frac{\zeta(2p)}{{\pi}^{2p}}=A_{2p}
\end{equation}

\begin{equation}
 \frac{\zeta(4n-1-2p)}{{\pi}^{4n-1-2p}}=B_{4n-1-2p}
\end{equation}

\begin{remark}
On sait que $A_{2p}$  est rationnelle
\end{remark}
\par
Ces trois d\'efinitions dans (2.38) alors on \'ecrit:
\begin{equation}
 B_{4n-1}=\sum_{p=1}^{2n-2}\frac{A_{2p}}{2n-1}.B_{4n-1-2p}
\end{equation}\\

Au lieu d' \'etudier  la nature de $\zeta(4n-1-2p)$ directement, nous sugg\'erons d'\'etudier d'abord celles de $B_{4n-1-2p}$ et $B_{4n-1}$; $n=2,3,4,...$ \`a partir de la relation ci-dessus.
\\
\\En effet, on peut \'ecrire la relation:
\begin{equation}
 B_{4n-1}=\sum_{p=1}^{2n-2}\frac{A_{2p}}{2n-1}.B_{4n-1-2p}
\end{equation}\\ sous forme matricielle comme suite:
\begin{displaymath} 
   \begin{pmatrix} 
         {B_{4n-1}}\\{B_{4n-5}}\\{B_{4n-7}}\\{...}\\{B_3}  
   \end{pmatrix}=
   \begin{pmatrix} 
         {\frac{A_{2}}{2n-1}}&{\frac{A_{4}}{2n-1}}&{\frac{A_{6}}{2n-1}}&{...}&{\frac{A_{4n-4}}{2n-1}}\\
         {0}&{1}&{0}&{...}&{0}\\
         {0}&{0}&{1}&{...}&{0}\\
         {...}&{...}&{...}&{...}&{...}\\
         {0}&{0}&{0}&{...}&{1}
   \end{pmatrix}
   \begin{pmatrix} 
         {B_{4n-3}}\\{B_{4n-5}}\\{B_{4n-7}}\\{...}\\{B_3}  
   \end{pmatrix} 
\end{displaymath}\\
La matrice ci-dessus est une matrice triangulaire. On sait qu'une matrice triangulaire a pour d\'eterminant le produit des \'el\'ements diagonaux. Comme aucun de ces  \'el\'ements n'est nul, alors cette matrice est toujours inversible et les deux familles(colonnes) ci-dessus ont la m\^eme nature; d\'ependante et d\'ependante(ou  ind\'ependante et ind\'ependante).  Nous parlons ici bien sure de la d\'ependance et de l'ind\'ependance lin\'eaire.
%\section{Valeurs num\'eriques de $\zeta(4n-1)$}
% \includepdf[pages=1-3]{CalculsDepuisZeta}
\section{Conclusion} 
Nous avons trouv\'e la fonctionnelle:
\begin{equation}
\zeta(4n-1)=\frac{1}{2n-1}\sum_{p=1}^{2n-2}\zeta(2p)\zeta(4n-1-2p)
\end{equation}\\$n=2, 3, 4, 5, 6, ...$\\
\\Elle est jusque-l\`a,  la formule condens\'ee de zeta aux impairs positifs. A partir d'elle, nous avions introduit une nouvelle demarche pour arriver \`a \'etudier la nature de zeta \`a ces entiers. La relation (3.5) nous montre que les constantes $ B_3, B_5,..B_{4n-3}$  cachent une propri\'et\'e pour $ \zeta(3),\zeta(5),.,\zeta(4n-3)$; $n=2,3,...$ 
\section{R\'ef\'erences}
Nous avons pr\'ef\'er\'e les d\'emonstrations que la citation des r\'ef\'erences.
\end{document}